\documentclass[11pt, reqno]{amsart}
\usepackage{epsf} 
\usepackage{amsmath, amsthm, amssymb, amsfonts, verbatim, color, float}
\usepackage[pdftex]{graphicx}
\usepackage[bookmarks]{hyperref}

\reversemarginpar            

\usepackage{caption}
\usepackage{subcaption}
\usepackage{color}

\newcommand{\R}{\mathbb{R}}

\def\R{{{\mathbb R}}}

\def\bv{{\bf{v} }}
\def\bx{{\bf{x} }}

\newcommand{\norm}[1]{\left\Vert#1\right\Vert}
\newcommand{\abs}[1]{\left\vert#1\right\vert}
\newcommand{\ip}[1]{\langle#1\rangle}
\def\be{\begin{equation}}
\def\ee{\end{equation}}
\def\eps{{\epsilon}}

\newtheorem{theorem}{Theorem}

\newtheorem{remark}[theorem]{Remark}

\newtheorem{lemma}[theorem]{Lemma}

\newtheorem*{acknowledgment}{Acknowledgments}

\numberwithin{equation}{section}
\numberwithin{theorem}{section}

\title[] { Blowing Up Solutions to the Zakharov System for Langmuir Waves}
\author{Y. ~Cher}
\address{Department of Mathematics, University of Toronto}

\author{M. ~Czubak}
\address{Department of Mathematical Sciences, Binghamton University (SUNY)}
\author{C.~Sulem}
\address{Department of Mathematics, University of Toronto}


\begin{document}

\begin{abstract}
Langmuir waves take place in a quasi-neutral  plasma  and are modeled
by the Zakharov system.  The phenomenon of collapse, described by blowing up solutions plays a 
 central role in their dynamics. 
We present in this article  a review of the main mathematical properties of blowing up solutions. They include conditions for 
blowup in finite or infinite time, description of  
self-similar singular solutions and lower bounds for the rate of blowup of certain norms associated to the solutions.

\end{abstract}

\maketitle

\section{Introduction}
Langmuir waves take place in a non-magnetized or weakly magnetized plasma and are described
by the Zakharov system \cite{Z72}
  \begin{align}
& i\partial_t {\bf E} - \alpha\nabla \times (\nabla \times {\bf E})  +  
\nabla(\nabla \cdot {\bf E}) =n {\bf E}, \label{Z1}  \\
&\partial_{tt} n - \Delta n = \Delta |{\bf E}|^2  \label{Z2}
\end{align}
Equation \eqref{Z2} originates from  the hydrodynamic system
\begin{align}
& n_t+\mbox{\boldmath $\nabla$}\cdot{\bf v}=0, \label{Z3}  \\
&{\bf v}_t+\mbox{\boldmath $\nabla$}n
= - \mbox{\boldmath $\nabla$}|{\bf E}|^2, \label{Z4}
\end{align}
governing ion sound waves.
${\bf E}({\bf x} ,t)$ is the  complex envelope of the electric field oscillations.
   $n({\bf x},t)$ denotes the  fluctuations of density of ions and ${\bf v}( {\bf x},t)$
 their velocity,  with ${\bf x} \in \R^d$, in dimension $d=$ 2 or 3.  
The parameter $\alpha$ in \eqref{Z1} is defined as the square ratio of the light speed and 
the electron Fermi velocity and is usually large.
A simplified  system of equations is obtained in the electrostatic limit ($\alpha \to \infty$) expanding 
the electric  field in the form  ${\bf E}  = \nabla \psi + \frac{1}{\alpha} {\bf E_1} +....$. Substituting this expansion
 in \eqref{Z1} and
taking the divergence of the equation gives a  system describing the interaction of the electrostatic
potential with the plasma density 
 \cite{Z72},\cite{ZMS75} 
 \begin{align}
& \Delta (i\partial_t \psi + \Delta \psi) =\nabla \cdot (n \nabla \psi), 
\label{zakhpot1}\\
&  \partial_{tt}n - \Delta n =\Delta (|\nabla \psi|^2).
\label{zakhpot2}
\end{align}
A further simplification leads to 
\begin{align}
i\psi_t+\Delta\psi = n\psi,  \label{z1} \\
 n_{tt}-\Delta n=\Delta |\psi|^2  \label{z2}.
\end{align} 
with
\begin{align}
 n_t+\mbox{\boldmath $\nabla$}\cdot{\bf v}=0,  \label{z3}\\
{\bf v}_t+\mbox{\boldmath $\nabla$}n
= - \mbox{\boldmath $\nabla$}|\psi|^2, \label{z4}
\end{align} 
 usually called the scalar Zakharov system.
Introducing the hydrodynamic potential $U$ such that ${\bf v} = - \nabla U$, eq.\eqref{Z4} becomes
\begin{equation}
\partial_t U = n + |{\bf E}|^2.
\end{equation}

Heuristic derivations of the Zakharov system can be found in \cite{R97}, \cite{B98}.
Viewing the plasma as a two interpenetrating fluids (electrons and ions), the Zakharov system \eqref{Z1}-\eqref{Z3}  
can be obtained using a multiple-scale modulation analysis \cite{SS99}. 
A rigorous derivation of the scalar model is given in \cite{T07} using techniques of geometric
optics and semi-classical calculus. 

Invariance properties of the system by simple transformations lead to several conserved quantities. In particular,
if $({\bf E}, n)$ is a smooth solution of  (\ref{Z1})--(\ref{Z4}), 
the wave energy 
$ N = |{\bf E}|_{L^2}^2 $ and the Hamiltonian
\begin{equation}
H = \alpha \,|\nabla \times {\bf E}|_{L^2}^2 + |\nabla \cdot {\bf E}|_{L^2}^2
+ \frac{1}{2} \,|n|_{L^2}^2 +  \frac{1}{2} \,|\nabla U|_{L^2}^2
+ \int n |{\bf E}|^2 d\,{\bf x}
\end{equation}
are conserved. Other invariants are the linear
 and angular momenta
\begin {equation}
{\bf P} = \int \Bigl(\frac{i}{2} \sum_j (E_j \nabla E_j^* - E_j^* \nabla 
E_j)+ n {\bf v}\Bigr)\,d{\bf x}
\end{equation}
and  
\begin{equation}
{\bf M} = \int (i {\bf E} \times {\bf E}^* + {\bf x} \times {\bf P})\,d{\bf x}.
\end{equation} 
Modulational instability leads to the formation of regions where the 
density of the plasma is very low. In these  regions referred to as cavities,  high-frequency oscillations of the 
electric field are trapped. Their nonlinear evolution gives rise  to the collapse of the cavities and a strong amplification of
  the amplitude of the oscillations of the electric field. Heuristic arguments
and numerical simulations show that, for large enough initial conditions,  solutions blow-up
in a finite time  both in two and three dimensions (see \cite{SS99} for a review).

In this article, we present an overview of mathematical results and open questions 
concerning blowing up solutions for the scalar  Zakharov model. We also discuss the extension of 
 some of the features  of blowup to the Vectorial Zakharov system for which 
very few rigorous are known apart from local wellposedness  and global 
wellposedness  under the assumption of  small enough initial conditions.    

\section{The scalar Zakharov system}

We consider the scalar Zakharov system \eqref{z1}-\eqref{z2} 
with initial conditions
\begin{equation} \label{IC} 
\psi(\bx,0)= \psi_0(\bx), \;  n(\bx,0)= n_0(\bx), \; 
n_t(\bx,0)= n_1(\bx).
\end{equation}
The conserved quantities are:\\
\noindent
the wave energy 
\begin{equation}   
N=|\psi|^2_{L^2}, 
\end{equation}
the linear momentum
\begin{equation}
{\bf P} = \int  \Big(\frac{i}{2} (\psi \nabla \psi^* - \psi^* \nabla
 \psi) + n{\bf v}\Big ) d{\bf x},
\end{equation}
the angular momentum
\begin{equation}
{\bf M} = \int {\bf x} \times {\bf P} \,d{\bf x},
\end{equation}
and the Hamiltonian 
\begin{equation}
H=\int\Big( |\mbox{\boldmath $\nabla$}\psi|^2+n|\psi|^2+
    \frac{1}{2}|{\bf v}|^2+\frac{1}{2}n^2 \Big)d{\bf x}. \label{HZAK}
\end{equation}
There is a large literature devoted to the local and global wellposedness
of the initial value problem.  Earlier works concern  smooth solutions, in particular  solutions with finite  energy (Hamiltonian)
(\cite{SS79}, \cite{AA84}, \cite{SW86}, \cite{OT92}, 
\cite{GM94}). \cite {BC96} . More recently, there has been an interest in solutions with lower 
regularity assumptions and in particular in  solutions with infinite energy  \cite{GTV97},  \cite{T00}, \cite{P01}, \cite{CHT08},
\cite{BHHT09}, \cite{BH11}.
Associated to the long time existence theory are the important questions  of  scattering theory,  existence of wave operators \cite{GV06},  and  precise decay of solutions for large time.  In particular, in three dimensions, it is proved in \cite{HPS13} that,
  if the initial conditions are small and localized, then $\sup_x|\psi(t) | \le C |t|^{-7/6-} $, $\sup_x |n(t)| \le C|t|^{-1}$,
   and the solution $(\psi,n)$ scatters to a  solution to the associated linear problem as $ |t| \to \infty$. Here the notation $7/6-$ means
   $7/6-\varepsilon$, for any $\varepsilon>0$.

We denote by $H^k(\R^d)$ the Sobolev space of functions $f$ such that $f$ and  its derivatives of  order $p$, $|p|\le k$, 
are bounded in the $L^2$-space. It is also convenient to define the product space 
\be\label{Hk}
H_k =H^k(\R^d) \times H^{k-1}(\R^d) \times H^{k-2}(\R^d). 
\ee
The energy space corresponds to $H_1$.
\subsection{Blowup in  finite or infinite time}

A central tool in the theory of blowup for the Nonlinear Schr\"odinger (NLS) equation
\begin{equation} \label{NLS}
i\partial_t\psi + \Delta \psi + |\psi|^{2\sigma} \psi=0, \;\; \psi({\bf x},t) =\psi_0({\bf x})
\end{equation}
 is  the variance identity 
\begin{equation*}
\frac{d^2}{dt^2} \int |{\bf x}|^2 |\psi|^2 d{\bf x}=8H_{NLS} - 4\frac{d\sigma -2}{\sigma + 1 } 
           \int |\psi|^{2 \sigma + 2} \,d {\bf x}, \label{Variance}
\end{equation*}
where $H_{NLS} = \int  \Big ( |\nabla \psi|^2 - \frac{1}{\sigma +1}
  |\psi|^{2\sigma +2}  \Big ) d{\bf x}$ is the NLS Hamiltonian. Under the assumption that the initial condition
  $\psi_0$ is in $H^1(\R^d)$, has  finite variance and $H_{NLS}(\psi_0) <0$, the solution of \eqref{NLS} blows up in a finite time
  if   $\sigma  d\ge 2$.
For  the Zakharov system, 
the usual variance
$\int |\bx|^2 |\psi|^2 d\bx$, can be  replaced by the quantity
\begin{equation}
  {\mathcal V}(t) \  = \  \frac{1}{4} \int |\bx|^2 |\psi|^2 d\bx + 
 \int_0^t \int ({\bf x} \cdot {\bf v} ) n\,d\bx\,dt \label{zhkvar}
\end{equation}
which is well-defined for functions in the space 
\begin{equation}
\Sigma' = \{ (\psi, n,{\bv}) \in H_1, 
\int \Big ( |\bx|^2 |\psi|^2 + |\bx| (|n|^2
+ |\bv|^2 )\Big)  d{\bf x} <  \infty \}. 
\end{equation}
The function ${\mathcal V}(t)$
 and satisfies
\begin{equation}
\frac{d^2{\mathcal V}}{dt^2}(t) =d H - (d-2) |\nabla \psi|_{L^2}^2   - 
(d-1)  |{\bf v}|_{L^2}^2 ,
\label{v22}
\end{equation} 
where $H$ is the Hamiltonian defined in \eqref{HZAK}.
 In dimension $d\ge 2$, one has
 $\displaystyle{\frac{d^2 {\mathcal V}}{dt^2} < 0}$ if the initial conditions are such that  $H< 0$. However,
one cannot conclude  on existence of blowup
solutions  because,  unlike the NLS case,  ${\mathcal V}$ does not have a fixed sign. In particular, in dimension 2,  
Merle \cite{M96a} proved that it tends to $ -\infty$ as the singularity is approached.
Nevertheless, one can get partial results, under the assumption of radial symmetry.
Indeed in this case, there is a useful result referred to as  the Strauss Lemma \cite {S77}
which gives an upper bound of  the sup norm of a function in term of  its  $H^1$-norm far from the origin. Namely, if $f$ is a radially symmetric function in $H^1(\R^d)$ with $d \geq 2$,
then, for any $R > 0$,
\begin{equation}
|f|_{L^\infty(|\bx|>R)}^2 \leq C R^{-d+1} |\nabla f|_{L^2(|\bx|>R)}
|f|_{L^2(|\bx|>R)}.
\label{lem}
\end{equation}
The radial assumption has  been useful in other contexts such as  in existence and scattering theory, where
it  allows a larger  range of parameters for linear estimates of  Strichartz type \cite{GN13}.

 The method consists in modifying the quadratic  weight $|\bx|^2$ 
in $\mathcal{V}$ defined in \eqref{zhkvar}
by a smooth function $p(\bx)$ that behaves like  $|\bx|^2$ near 
the origin and like $|\bx|$ at infinity, and by considering the time derivative
$\displaystyle{y(t) =  -\frac{d{\mathcal U}}{dt} }$ of the modified variance
\begin{equation}
{\mathcal U}(t)  =  \frac{1}{2} \int p(x) |\psi|^2 ~d\bx +
\int_0^t \int (\nabla p  \cdot {\bf v})\,n\,d\bx\,dt,
\end{equation}
as in the case of solutions of the NLS equation with infinite variance. 
However, the modification of the weight induces additional terms in the  time evolution of 
the function $y (t) $ that need to be estimated.
For this purpose, one uses a  sequence of rescaled weights $p_m( |\bx|) = m^2 p(\bx/m)$, and  proves that the
additional contributions are controlled for $m$ sufficiently large.
More precisely,  one proves that
\begin{equation}
y_m (t) = - \Im \int ((\nabla p_m \cdot \nabla \psi) \psi^* - 
 (\nabla p_m \cdot  {\bf v}) n)\,d{\bf x} 
\end{equation}
satisfies, for $m$  sufficiently large
\begin{equation}
y_m (t) \ge \frac{d}{2} |H| t.
\end{equation}
On the other hand, the function $y_m(t)$ is controlled by the norm of the solution in the energy space, namely
\begin{equation}
|y_m(t) | \le C (|\psi_0|_{L^2}^2 + |\nabla\psi|_{L^2}^2 +
|{\bf v}|_{L^2}^2 + |n|_{L^2}^2).
\end{equation}
We have the  following result  proved by Merle  in \cite{M96b}:
\begin{theorem} \label{thmvarz} \  
Consider the Zakharov system \eqref{z1}, \eqref{z3}, \eqref{z4}  in dimension $d=2$ or $d=3$  with initial conditions
in the space  $\Sigma'$.
Assume that there exists a  smooth solution $(\psi,n,{\bf v} )$  during an
interval of time $[0,t_0]$. In particular, its mass and  Hamiltonian  are conserved and its variance ${\mathcal U}$ is
well-defined. Assume in addition that
the solution  is radially symmetric and its Hamiltonian $H<0$. Then,
either
$|\psi|_{H^1} + |n|_{L^2} +  |{\bf v}|_{L^2}
 \rightarrow \infty $ as $t \rightarrow t^*$ with $t^*$ finite,
or $(\psi,n,{\bf v})$ exists for all time and 
 $|\psi|_{H^1} + |n|_{L^2} +  |{\bf v}|_{L^2} 
\rightarrow \infty $ as $t \rightarrow \infty$.
\end{theorem}
\begin{remark}
An open question is the extension of this analysis to solutions that are not radially symmetric, and furthermore,
to  solutions to the full vector Zakharov system.
Based on numerical observations, it is believed that blow up does  indeed occur in a finite time
 for general initial conditions with negative Hamiltonian.

\end{remark}
\subsection{Self-similar blowing up solutions}

\subsubsection{Dimension $d=2$}
Unlike for the NLS equation, there is no conformal mapping for the two-dimensional  Zakharov system. Nevertheless one can construct
 exact self-similar blowing up solutions that have the form
in the form
\cite{ZS81}
\begin{align}
&\psi( x,t)= \frac{1}{a (t^*-t)} P\Big(\frac{|x|}{a(t^* -t)}
\Big)
e^{i\left (\theta + \frac{1}{a^2(t^*-t)}-\frac{|x|^2}{4(t^*-t)}\right )}, \label{selfsimzak1}\\
&n(x,t) = \frac{1}{a^2(t^*-t)^2} N\Big(\frac{|x|}{a (t^*-t)}\Big),
\label{selfsimzak2}
\end{align}
where $(P,N$)  are real functions satisfying the system of ODEs 
\begin{align}
&\Delta P-P-NP =0, \label{Ppeq2} \\
&a^2(\eta^2 N_{\eta \eta}+6\eta N_{\eta}+6N)- \Delta N=\Delta P^2,
\label{Nneq2}
\end{align}
with  $\eta$ being  the rescaled independent variable and  $a>0$  a free parameter. 
Glangetas and Merle have rigorously  studied   the system    \eqref{Ppeq2}-\eqref{Nneq2} in  \cite{GM94}. We summarize below the most important properties. 
When $a=0$, $N= -P^2$ and  eq. \eqref{Ppeq2}  becomes
 \begin{equation}
 \Delta P-P+P^3 =0. \label{Req}
 \end{equation} 
 It is known that  Eq. \eqref{Req}  has an infinite number of radially symmetric solutions that decay exponentially at infinity, only one of them, denoted $R$ and called the NLS
the ground state,  (also known as the Townes soliton) 
is strictly positive and monotone decreasing (see for example \cite{BL83}). It  plays a central role in the study of NLS equations.
  
If the coefficient $a$  in \eqref{Nneq2} is sufficiently small, there exists a  solution $(P_a, N_a)$ in $H^1\times L^2$  with $P_a >0$.
This solution is in fact $C^\infty$ and  its derivatives of order $k$ satisfy the decay properties
\begin{equation*}
|P^{(k)}(\eta)|\le  c_k e^{-\delta\eta} , \; \; \quad N^{(k)}(\eta)| \le \frac{c_k}{|\eta|^{k+3}}
\end{equation*}
for large $\eta$.
When the parameter $a$ is small,  the solution $(P_a,N_a)$ is constructed by a continuation method from the solution $(R, -R^2)$ corresponding to $a=0$.  

Furthermore, for any value $c$ strictly larger than the $L^2$-norm of the NLS ground state $R$, 
there exists $a_c$ such that for any $a<a_c$, there is a unique  solution $(P_a,N_a)$  in $H^1\times L^2$
with $P_a>0$ and $|P_a|_{L^2} < c$. 

Numerical simulations show that for a large class of radially symmetric  initial conditions having a strictly negative Hamiltonian, the solutions display a self-similar collapse as $t\to t^*$  as described by \eqref{selfsimzak1}-\eqref{selfsimzak2} \cite{BPP90}, \cite{LPSSW92}. The coefficient $a$ in the equation 
for the limiting profiles $(P,N)$ depends on the initial conditions. 
 When considering a sequence of
initial conditions with an initial $L^2$-norm of $\psi_0$ decreasing to 
$|R|^2_{L^2}$,  where $R$ is the NLS ground state, it was  observed that the computed value of the coefficient $a$ 
tends to zero. In this limit, the self-similar
profile becomes (strongly) subsonic and 
tends to the NLS ground state  $R$. This limit is
delicate, since  solutions of the scalar Zakharov equation  with  critical norm $|\psi_0|_{L^2}^2 = |R|_{L^2}^2$ remain smooth
for all time \cite{GM94}.
Indeed, unlike the NLS equation, there are no minimal blowing up solutions to the 2d Zakharov system.
For initial conditions in the energy space such that $|\psi_0|_{L^2}\le |R|_{L^2} $  solutions remain in the energy space for all times.  The case  $|\psi_0|_{L^2} < |R|_{L^2} $ 
is straightforward and follows the NLS analysis   \cite{AA84}, \cite{W83}.  When  $|\psi_0|_{L^2} = |R|_{L^2} $, the global wellposedness property is very specific to the Zakharov system.

Finally, when numerical simulations  are performed with anisotropic initial conditions with negative Hamiltonian, it was
observed that the solutions become isotropic near collapse with the same limiting profiles as those obtained with isotropic initial conditions \cite{LPSSW92}.

\subsubsection{Dimension $d=3$}
In three dimensions, there are no known explicit blowing up solutions.
Self-similar solutions exist only asymptotically close to collapse and have the universal form
\cite{BZS75}, \cite{ZS81}

\begin{align}
&\psi(x,t) \sim  \frac{1}{(t^*-t)}P\Big(\frac{|x|}{\sqrt{3}(t^*-t)^{2/3}}\Big)
e^{i(t^*-t)^{-1/3}}, \label{r.ss}\\ 
&n(x,t)  \sim \frac{1}{3(t^*-t)^{4/3}}
N\Big(\frac{|x|}{\sqrt{3}(t^*-t)^{2/3}}\Big),\label{n.ss} 
\end{align}
where $P(\eta)$ and $N(\eta)$ are radially symmetric scalar functions satisfying 
the coupled system of ODEs
\begin{align}
&\Delta P-P-NP=0,  \label{P3eq} \\
&\frac{2}{9}(2\eta^2 N_{\eta\eta} + 13 \eta N_{\eta}
+ 14 N)=\Delta {P}^2  \label{N3eq}.
\end{align}
This type of blowup is referred to as supersonic collapse, because, when substituting the expressions \eqref{r.ss}-\eqref{n.ss} into
the Zakharov system, the pressure  term
 $\Delta n$ is of lower order than $\partial_{tt}n $.

Note that, unlike the 2d case,  there is no free parameter in the system. 
As discussed in  \cite{ZS81} and proved more recently in \cite{M01a},
there exists an infinite number of solutions $(P_k,N_k)$ to \eqref{P3eq}-\eqref{N3eq} 
such that, for all $k$,   
\begin{equation}
0<P_k(\eta) < P_{k+1} (\eta) \;  \; \rm{and} \; \;  N_k(\eta) <0.
\end{equation}
The profiles $P_k$ decay exponentially $|P_k (\eta)| \le  C_k e^{-\delta \eta} $ for $\eta$ large, while the $|N_k| \le
\frac{C_k}{1 +\eta^2}$ decay algebraically.

The values of $P$ and $N$ at the origin satisfy the relation
$N_k(0) = \frac{9 P_k(0)^2}{14-9P_k^2(0)}$. We have also $P'(0)=N'(0)=0$ due to the radial symmetry.
The pair  $(P_k,N_k) $ is thus characterized by  the value 
$P_k(0)$.  It is proved in \cite{M01a} that  there exists a sequence $\alpha_k =   \frac{1}{3} \sqrt{2k(4k+3)}>0$ such that the values 
$P_k(0)$  are ordered as 
$\alpha_k< P_k(0) < \alpha_{k+1} $. 
It is of interest to see how the values $\alpha_k$ arise in the analysis. They  appear when
 one writes the Taylor series  expansion of $P$ and $N$ near the origin:
  \begin{equation}
  P(\eta) = \sum\limits_{i=0}^\infty a_i \eta^{2i} \ ; \ N(\eta) = \sum\limits_{i=0}^\infty b_i \eta^{2i}.
  \end{equation} 
  The series have only even powers because $P$ and $N$ are radially symmetric.
 From the substitution of  the Taylor series into the system \eqref{P3eq}-\eqref{N3eq}, one  gets  two  relations between the coefficients $a_i, b_i$.  The values $\alpha_i$ appear  when solving the equation for the coefficient $a_i$:
 \begin{equation}
 (\alpha_i^2 - a_0^2) a_i = F( a_0, ...., a_{i-1}, b_0, ....b_{i-1})
 \end{equation}
 when solving for the coefficients.
  In order to have well defined coefficients and an analytic solution,  $P(0)$ which identifies to $a_0$ should be different
 from  the $\alpha_i$. In \cite{M01a}, it is proved that there is at least one solution $P_k$ with initial value $P_k(0) \in (\alpha_k,\alpha_{k+1})$
 which is strictly positive and decays to $0$ at infinity.  
 Numerically, we found (at least for the first few that we computed) that
 there is only one.
Figs.1a and 1b  show the first four pairs of solutions computed numerically
by a shooting method (with the shooting parameter being $P_k(0)$. The values $P_k(0)$ (for $k=1,..,4$) are:
\begin{align}\label{eqn:p_0-vals}
P_1(0) \approx 1.38, \quad P_2(0) \approx 2.43 ,\nonumber \\
P_3(0) \approx 3.42,  \quad P_4(0) \approx 4.40.
\end{align}
\begin{figure}[tbp]
\begin{center}
\begin{subfigure}{0.5\textwidth}
  \begin{center}
 \includegraphics[scale=.29]{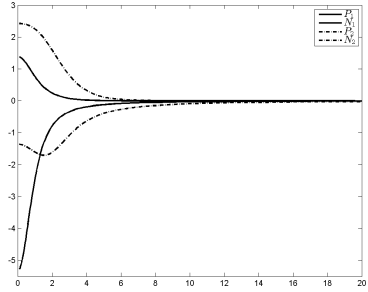}
  \end{center}
\end{subfigure}
\begin{subfigure}{.5\textwidth}
\begin{center}
 \includegraphics[scale=.29]{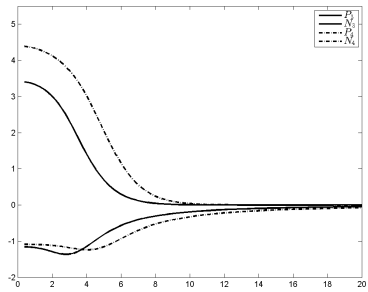}
  \end{center}
\end{subfigure}
\end{center}
\caption{Solutions $(P_k, N_k)$ of \eqref{P3eq}-\eqref{N3eq} for $k = 1,...,4$. Top: Solid line -- $(P_1, N_1)$, dashed line -- $(P_2, N_2)$ corresponding to initial values $P_1(0)$ and $P_2(0)$ in \eqref{eqn:p_0-vals} respectively. Bottom: Solid line -- $(P_3, N_3)$, dashed line -- $(P_4, N_4)$ corresponding to initial values $P_3(0)$ and $P_4(0)$ respectively.}
\label{fig:ZakhSelfSim}
\end{figure}


Like in the 2d case, 
the dynamical stability of the  asymptotically
self-similar solutions to the 3d Zakharov system 
for both radially symmetric and anisotropic 
initial conditions  was studied numerically in \cite{LPSSW92}.
It was observed that for a large class of data, blowup
solutions   asymptotically  display a self-similar collapse  described by  the above solutions.
The profiles identify to the first mode $(P_1,N_1)$ solution of \eqref{P3eq}-\eqref{N3eq}
that has the lowest value at the origin, and for which $N_1$ is monotone increasing.

\begin{remark} 
There is no rigorous proof of  dynamic stability of the (2d) self similar  or (3d)  asymptotically self similar  
solutions even for well prepared initial conditions  (with or without radial symmetry)
chosen close to the profiles $(P,N)$ solutions of  the ODE systems 
\eqref{Ppeq2}-\eqref{Nneq2} or \eqref{P3eq}-\eqref{N3eq}.
\end{remark}

\subsection{Lower bounds for rate of blowup}

\subsubsection{Scale invariance, criticality and local wellposedness}

 An important aspect  in the analysis of dispersive equations is the notion of criticality. It is closely 
 related to the invariance properties  of the equation.
 For example, the  NLS equation  \eqref{NLS} is invariant under the scaling transformation
 $ \psi({\bf x}, t) \rightarrow \psi_\lambda({\bf x},t) =  \lambda^{1/\sigma} \psi(\lambda {\bf x}, \lambda^2 t)$.
 It is said to be  $\dot H^s$-critical if the (homogeneous) $H^s$-norm  is unchanged under the 
 above scaling transformation.  The  corresponding critical Sobolev exponent for NLS is  thus $s_c= d/2 -1/\sigma$. The notion of criticality
is not straightforward for the  Zakharov system because  the Schr\"odinger equation and the wave equation have 
different scale invariances.  
In \cite{GTV97},  criticality is defined
by considering the scaling
\begin{equation}
\psi \to \psi_\lambda = \lambda^{3/2} \psi(\lambda \bx, \lambda^2t), \quad
 n \to n_\lambda =\lambda^2 n(\lambda \bx, \lambda^2t)
\end{equation}
that would leave the Zakharov system invariant in the absence of the 
term $\Delta n$. This is indeed  the
relevant scaling to study  blowing up solutions of the three dimensional 
Zakharov system as we have seen in the previous section. 

In relation to the initial value problem, the Sobolev space with critical exponent often corresponds to 
the space with minimal regularity in which the problem is locally
well-posed.  For the Zakharov system, the critical values for the initial value  problem  in 
$H^k \times H^l \times H^{l-1}$ are $ k= \frac{d}{2} -\frac{3}{2}$ 
and $ l=\frac{d}{2}-2$.
Note that  $k-l=\frac{1}{2}$,  while one would have $k-l=1$ in the classical setting of the energy space  $H_1$.
We now summarize the wellposedness results from the works of   \cite{GTV97}, \cite{CHT08}, \cite{BHHT09}, \cite{BH11}.  For ill-posedness results in dimension one, see \cite{H07}.
\begin{theorem}
In dimension $d=1$,  the Zakharov system is localy wellposed
in $H^k \times H^l \times H^{l-1} $, provided that
$ -\frac{1}{2} < k-l \leq 1, \; \; \; 2k \geq l+\frac{1}{2} \geq 0$. Furthermore,
global wellposedness  holds in the largest space in which local wellposedness holds, that is $L^2 \times H^{-1/2} \times H^{-3/2} $.

In dimension 2,
it is locally well-posed in the critical space  $L^2 \times H^{-1/2} \times H^{-3/2}$,  and in dimension 3,
it is locally well-posed in $H^{\varepsilon} \times H^{-1/2+\varepsilon} \times H^{-3/2+\varepsilon}$ which is also, up
to arbitrarily small $\varepsilon$ the critical space.

Finally, in dimension $d\geq 4$, the whole range of subcritical values
$k > \frac{d}{2}-\frac{3}{2}$ and $l> \frac{d}{2}-2$
 is covered by the theorem as long as $l\leq k\leq l+1$ and $2k-l-1>\frac d2 -2$.
\end{theorem}

\subsubsection{Finite energy solutions: the two-dimensional case.}

The next theorem is due to Merle \cite{M96a}. It concerns solutions of the 2d scalar Zakharov system with initial conditions
$(\psi_0,n_0,n_1)$ in the 
energy space  $H_1$, thus having  a  finite Hamiltonian. Assume  that there exists  a finite  time $t^*$  such that,
$$
|\nabla \psi(t)|_{L^2}  + |n(t)|_{L^2} + |{\bf{v}}(t)|_{L^2} \to \infty \;\; {\rm as} \;\; t\to t^*.
$$
The question is to determine at what rate these norms become infinite as $t$ approaches $t^*$.
\begin{theorem}  
Assume that the solution $(\psi,n)$  to the  2D scalar Zakharov system blows up  in the energy space  $H_1$ at a finite time $t^*$.  Then
there exist constants $c_1 >0$ and $c_2 >0$ 
depending only on $|\psi_0|_{L^2}$ such that for $t $ close  to $t^*$,
\begin{align}
&|\nabla \psi(t)|_{L^2}  \geq  \frac{c_1}{t^*-t}, \label{blpsi} \\
&|n(t)|_{L^2}  \geq \frac{c_2}{t^*-t}.
\label{bln}
\end{align}
More precisely,
the constants $c_1$ and $c_2$ 
scale like $(|\psi_0|_{L^2}^2 - |R|_{L^2}^2)^{-1/2}$ 
where $R$ is the NLS ground state.
\end{theorem}

\begin{remark}
This rate is  optimal in the sense that 
 the self-similar solutions \eqref{selfsimzak1}-\eqref{selfsimzak2} satisfy
\begin{equation}
|\nabla \psi(t)|_{L^{2}} = \frac{1}{a(t^*-t)} |\nabla P|_{L^2} , \qquad |n(t)|_{L^2} = \frac{1}{a(t^*-t)} |N|_{L^2}.
\end{equation}
and  thus  blow up exactly at the rate stated in theorem. Notice also that the theorem provides the blowup rate for the two  quantities 
$ |\nabla \psi(t)|_{L^{2}} $ and $|n(t)|_{L^2}$ separately but does not give information on  $|{\bf{v}}(t)|_{L^2} $.
\end{remark}

The derivation of this result is based on  scaling properties and conservation of the Hamiltonian. 
 One defines the   rescaled functions ${\widetilde \psi}({\bf x},s) ,
{\widetilde  n}({\bf x}, s),{\widetilde {\bf v}}({\bf x},s) $ (where $t$ is seen as a parameter)
\begin{align}
&{\widetilde \psi}({\bf x},s) = \frac{1}{\lambda(t)}\psi
\Bigl(\frac{{\bf x}}{\lambda(t)}, 
t +\frac{s}{\lambda(t)}\Bigr), \\
& {\widetilde  n}({\bf x}, s) = \frac{1}{\lambda^2(t)} 
n\Bigl(\frac{{\bf x}}{\lambda(t)}, 
t +\frac{s}{\lambda(t)}\Bigr), \\
& {\widetilde {\bf v}}({\bf x},s) = \frac{1}{\lambda^2(t)}{\bf v}
\Bigl(\frac{{\bf x}}{\lambda(t)}, t +\frac{s}{\lambda(t)}\Bigr).
\end{align}
where the scaling factor 
\begin{equation}
\lambda (t) = \int \Big (|\nabla \psi|^2 + \frac{1}{2} n^2 + \frac{1}{2}
|{\bf v}|^2 \Big ) d{\bf x},
\end{equation} 
is associated to the energy norm. Notice that the scaling of the time variable corresponds to the wave equation rather than to the
Schr\"odinger equation.
At $s=0$, 
\begin{equation} 
\int \Big (|\nabla \widetilde \psi(0)|^2 + \frac{1}{2} \widetilde n(0)^2 + \frac{1}{2}
|{\bf \widetilde v}(0)|^2 \Big ) d{\bf x} =1.
\end{equation}
Under the hypothesis of the theorem,
$\lambda(t)\to \infty$  as $t$ approaches $t^*$.

The analysis consists in establishing bounds for the individual quantities
$|\nabla \widetilde \psi(0)|$, $\widetilde n(0)$, $\widetilde {\bf v}(0)$
and estimates of $ \widetilde \psi(s)$, $\widetilde n(s)$ and $\widetilde {\bf v} (s)$ as $t\to t^*$.
It  uses delicate compactness arguments allowing the identification of limiting quantities as $t$ goes to $t^*$.
This approach, initiated in \cite{M96a} and  now known as profile decomposition, has led to many breakthroughs  in various  fields of dispersive PDEs.

\subsubsection{Infinite energy solutions}
We present in this section another approach for the derivation of a lower bound for the rate of blowing up solutions.
It  is more general, but less precise than the one presented in the previous section. It  applies to the  problem in  dimensions
two or three, and to initial conditions that may  or may not have a finite  Hamiltonian.  The result below was established in 3D in \cite{CCS13}.  The 2D result follows the same line of proof.

\begin{theorem} \label{thmrate}
Let the  initial data $(\psi(0), n(0), n_t(0))$
 be  in    $ {\mathcal H_\ell} := H^{\ell+1/2}(\R^d)  \times H^{\ell}(\R^d) \times H^{\ell-1}(\R^d)$, with the condition
  $0\leq \ell \leq 1/2 $ if $d=2$ and $0\leq \ell \leq 1$ if $d=3$.
  Assume that the solution $(\psi,n,n_t)$ blows up in a finite time
$t^*$, that is, as $t$ approaches $t^*$, $\|\psi(t)\|_{H^{\ell+1/2}} +\|n(t)\|_{H^{\ell}} + \|n_t(t)\|_{H^{\ell-1}} \to \infty$.
Then,  the rate of blowup   of the 
 Sobolev norms satisfies the lower bound estimate
\begin{equation}
\|\psi(t)\|_{H^{\ell+1/2}} +\|n(t)\|_{H^{\ell}} + \|n_t(t)\|_{H^{\ell-1}}
> C(t^*-t)^{-\theta_\ell} \label{BUPrate3d}
\end{equation}
with $\theta_\ell = \frac{1}{4}(4- d +2 \ell)^-$ , $d=2$ or $d=3$.
\end{theorem}
In the above formula, the notation  $a^-$ means $(a-\varepsilon)$ for  arbitrarily small  $\varepsilon >0$.  
\begin{remark}
Unlike the method of the previous section, this approach provides  a lower bound for the sums of the norms 
of $\psi, n, n_t$ but not  for the norms  separately. 
\end{remark}

\begin{remark}
In 3D, the lower bound \eqref{BUPrate3d} is probably not  optimal.
 Indeed,
the homogeneous $\dot H^{\ell+1/2}$-norm of $\psi$ and the homogeneous $\dot H^{\ell}$-norm of $n$ in the expression
of the asymptotic solution \eqref{r.ss}-\eqref{n.ss} both blowup at  a faster rate, namely $\frac{1}{3}( 1+ 2 \ell)$. 
 Note that  these norms blowup at the same rate in 3D, showing that the  space $H^{\ell+1/2} \times  H^{\ell}$ is appropriate for the analysis.
\end{remark}

\begin{remark}
In 2D, the norms $\|\psi\|_{\dot H^k}$ , with $k=\ell +1$ and $\|n\|_{H^\ell}$ of the exact  self-similar solutions
\eqref{selfsimzak1}-\eqref{selfsimzak2} blow up at the same rate $(t^*-t)^{-(\ell+1)}$.  Merle's work \cite{M96a} gives the optimal rate when $k=1$. Theorem \ref{thmrate} predicts rates of blowup in the space ${\mathcal H_\ell}.$  It gives almost the optimal rate of blowup for $\psi$ when $\ell=0$, but it is off then by $\frac 12$ for $n$.
 
 \begin{remark}
In  3D, a particular result about blowup of a space-time norm $L^{q,r}_{x,t}$of $n$ is given in \cite{M01a} under the assumption that  blowup occurs  in the energy space $H_1$.


\end{remark}

\end{remark}

Assume that  the solution $u:=(\psi,n,n_t)$ exists during a finite time $|t|\le T$ in the space $\mathcal{H}_\ell$.
There are two elements in the proof
of the theorem above:  \\
(i) A local wellposedness estimate  for $u$ in the form of the one  obtained by Ginibre-Tsutsumi-Velo \cite{GTV97},
\be\label{apriori1}
\norm{u}_{X_T}\leq C\norm{u_{0}}_{\mathcal{H}_\ell }+CT^{\theta}\norm{u}^{2}_{X_T},\ \ \theta >0,
\ee
 where $\|\cdot \|_{X_T}$ is  a space-time norm that will be defined later. For the purpose of local wellposedness, it is not important to determine exactly the power  $\theta$, the only requirement being that it is away from 0. On the other hand, in the blowup analysis, 
a key element is  to maximize the power $\theta$,
 because it leads to a better estimate for  the  lower bound of the blowup rate.  
 We find an expression for $\theta$ that depends on $\ell$, the order of the norm under consideration and increases with $\ell$.

 (ii) A classical contradiction argument 
introduced in  \cite{W81}  for semilinear 
heat equations and used in \cite{CW90} for NLS equations that reverses the local wellposedness estimate into
a blowup rate estimate.

{\it Step 1:  Local wellposedness estimate.}
Rewrite the wave equation
\eqref{z2} as two reduced wave equations for  $$w^\pm =n\pm i \omega^{-1} \partial_t n,$$
 where  $\omega =(-\Delta)^{1/2}$.
   The Zakharov system then becomes 
 \begin{align}
i\partial_t \psi + \Delta \psi &=( w^++w^-)  \psi ,  \label{zakh1a}\\
(i\partial_t  \mp \omega) w^\pm &= \pm \omega ( |\psi|^2) . \label{zakh2a}
\end{align}
$(\psi, w^\pm)$ solve \eqref{zakh1a}-\eqref{zakh2a} with initial data $(\psi_{0}, w_{0}^\pm)=(\psi_{0}, n_{0}\pm i\omega^{-1}n_{1})$ if and only if $(\psi, n)$ solve \eqref{z1}-\eqref{z2} with initial data $(\psi_{0}, n_{0}, n_{1})$.

In the analysis, one slightly modifies the above system by replacing the operator $\omega=(-\Delta)^{1/2} $ by $\omega_1=(1-\Delta)^{1/2}$ to avoid divergence at low wavenumbers. This leads to an additional term in the wave equation of the form
 $\ip{\nabla}^{-1}\mathcal Re~w^\pm $, which  is linear with a gain in derivatives, thus it is easily controlled  (see \cite{CCS13}).

The solution of  \eqref{zakh1a}-\eqref{zakh2a}  is written in its Duhamel formulation. 
Since the solution is considered in a fixed interval $[-T,T]$, we introduce in addition a cut-off  $C^\infty$ function 
$\varphi(t) =1$ for $|t|\leq 1$,   $\varphi(t) =0$ for $|t|\ge 2$, $0\le \varphi(t)\le 1$, and define   $\varphi_T(t) =\varphi(t/T)$,
($T\le1$).
The initial value problem \eqref{zakh1a}-\eqref{zakh2a} on time interval $[-T,T]$  is equivalent to the system of integral equations
\begin{align}
\psi(t) &=\varphi_1(t) U(t) \psi_0-i\varphi_T(t)\int_0^t U(t-s)\varphi^2_{2T} (w^++w^-) \psi(s) ds,\label{d4psi}\\
 w^\pm(t)&= \varphi_1(t) W(t) w_{0}\pm i\varphi_T(t)\int_0^t W(t-s) \varphi^2_{2T}\omega ( |\psi|^2) ds,\label{d4n}
\end{align}
where $U(t) = e^{it \Delta},\  W(t) =e^{\mp it \sqrt{-\Delta}}$ are the free Schr\"odinger  and free reduced wave operators respectively.

The space $X_{T}$ in \eqref{apriori1}, in which the analysis performed, is a  product of weighted Sobolev spaces, with 
space-time  weights being the Fourier multipliers
associated to the linear Schr\"odinger and linear reduced wave equation \cite{B93}. Namely, $X_T= X_S^{\ell+\frac 12,b}\times X_{W_\pm}^{l,b},$ 
with the norms given by 
\begin{align*}
  \norm{\psi}_{X_S^{\ell+\frac 12,b}}&=\norm{\ip{\xi}^{\ell+\frac 12}\ip{\tau+\abs{\xi}^2}^{b}\hat\psi(\tau,\xi)}_{L^2_{\tau,\xi}},\\
\norm{w}_{X_{W_\pm}^{\ell,b}}&=\norm{\ip{\xi}^{\ell}\ip{\tau\pm\abs{\xi}}^{b}\hat w(\tau,\xi)}_{L^2_{\tau,\xi}},
\end{align*}
where $\ip{\xi} = (1 +|\xi|^2)^{1/2}$ and  $b>\frac 12$.

%
There are two distinct elements in the estimates, the linear estimates and the nonlinear ones.
\begin{lemma}(Linear estimates)  \cite[Lemma 2.1]{GTV97} \label{l:21}
Consider the general linear equation
\[
iu_t-\Phi(-i\nabla)u=F \ \mbox{on} \ [0,T]\times \R^d, \quad u(0)=u_0\in H^s,
\]
where $\Phi$ is a real valued function. Then for
$\frac 12-\eps<b\leq 1-\eps $, ($\eps >0$)
\begin{equation}\label{l:21e1} 
\|\varphi_T u \|_{X^{s,b}} \lesssim \|u_0\|_{H^s} +  T^{\eps}\| F\|_{X^{s,b-1+\eps}},
\end{equation}
where the norm  in $X^{s,b}$ is associated to the linear operator, namely
$\norm{\cdot}_{X^{s,b}}=\norm{\ip{\xi}^s\ip{\tau+\Phi(\xi)}^b\ \widehat \cdot(\tau,\xi)}_{L^2_{\tau,\xi}}$.\end{lemma}

An application of  \eqref{l:21e1}  to the solution $(\psi,w^\pm)$ gives
\be\label{dh1}
\norm{\psi}_{X_S^{\ell+\frac 12,b}}\lesssim \norm{\psi_0}_{H^{\ell+\frac 12}}+T^{\eps}\norm{\varphi^2_{2T}(w^++w^-) \psi}_{X_S^{\ell+\frac 12,b-1+\eps}},
\ee
and
\be\label{dh2}
 \norm{w^\pm}_{X_{W_\pm}^{\ell,b}}\lesssim \norm{w_{0}}_{H^{\ell}}+T^{\eps}\norm{\varphi^2_{2T}\omega\abs{\psi}^2}_{X_{W}^{\ell,b-1+\eps}}.
\ee

We now explain how one gets  an estimate for the nonlinear terms and produce higher powers of $T$.   The goal is to establish 
\begin{align} 
\norm{\varphi^2_{2T}w^+\psi}_{X_S^{\ell+\frac 12,b-1+\eps}}&\lesssim 
T^{\theta}\norm{\varphi_{2T} w^+}_{X_{W_+}^{\ell,b}}\norm{\varphi_{2T}\psi}_{X_S^{\ell+\frac 12,b}}, \label{m10}\\
\norm{  \varphi_{2T}^2\omega |\psi|^2 }_{X_{W_+}^{\ell,b-1+\eps}}
&\lesssim T^{\theta}
\norm{\varphi_{2T} \psi}_{X_S^{\ell+\frac 12,b}}^{2} , \label{m20}
\end{align}
with equivalent estimates for $w^-$.
 A classical argument is to consider the nonlinear terms on the Fourier side and use duality.  This reduces   \eqref{m10}-\eqref{m20} to showing the following inequalities
\begin{align}
\abs{N_1}\lesssim T^\theta\norm{v}_2\norm{v_1}_2\norm{v_2}_2\label{m1},\\
\abs{N_2}\lesssim T^\theta\norm{v}_2\norm{v_1}_2\norm{v_2}_2,\label{m2}
\end{align}
where
\begin{align}
N_1&=\int\frac{\hat v( \xi_{1}-\xi_2, \tau_{1}- \tau_2) \hat v_1(\xi_1,\tau_1)\hat v_2(\xi_2,\tau_2)~\ip{\xi_1}^k}
{\ip{\tau+\abs{\xi}}^b\ip{\tau_1+\abs{\xi_1}^2}^{c}\ip{\tau_2+\abs{\xi_2}^2}^{b}\ip{\xi_2}^k\ip{\xi}^l }  ~d\xi_1 d\xi_2 d\tau_1 d\tau_2 ,\\
N_2&=\int\frac{\hat v( \xi_{1}-\xi_2, \tau_{1}- \tau_2) \hat v_1(\xi_1,\tau_1)\hat v_2(\xi_2,\tau_2)~\abs{\xi}~\ip{\xi}^\ell }
{\ip{\tau+\abs{\xi}}^c\ip{\tau_1+\abs{\xi_1}^2}^{b}\ip{\tau_2+\abs{\xi_2}^2}^{b}\ip{\xi_1}^k\ip{\xi_2}^k }~d\xi_1 d\xi_2 d\tau_1 d\tau_2 .
\end{align}
Ginibre-Tsutsumi-Velo \cite{GTV97} showed the above estimates by a repeated application of an inequality obtained from Strichartz estimates and H\"older inequality in time (see \cite[Lemmas 3.1-3.4]{GTV97}).  Their analysis did not require an optimal power of $\theta$, but needed it to be just large enough, so the final power of $T$ was positive.  They find 
\[
(b+1-(\frac n2+1)b_0-\epsilon)\frac 1{2b}.
\]
 For the rate of blowup analysis, we seek the optimal power of  $\theta$, 
and obtain estimates \eqref{m1}-\eqref{m2} with
\begin{align}\label{temptheta}
\theta=b+1-(\frac n2+1-\ell)b_0-\epsilon.
\end{align}
To remove the time cut off from the right hand side of  \eqref{m10}-\eqref{m20}, we recall
\begin{lemma}\cite{GTV97}\label{l:22}
\[
\norm{\varphi_Tu}_{X^{s,b}}\leq CT^{-b+\frac 1q}\norm{u}_{X^{s,b}},
\]
where $s \in \R, b\geq 0, q\geq 2$ and $bq>1$.
\end{lemma} 
Applying this twice (since the nonlinearity is quadratic) with $q=2$ and combining estimates \eqref{dh1}-\eqref{dh2} with \eqref{m10}-\eqref{m20}  gives the final estimate
\begin{align}
\|\psi\|_{X_S^{\ell+\frac 12,b}}  + \|n\|_{X_{W_+}^{\ell,b}}+&\|n_{t}\|_{X_{W_+}^{\ell-1,b}} \le 
 C\big(\| \psi_0 \|_{H^{\ell+\frac 12}} +\|n_{0}  \|_{H^\ell}+  \|n_{1}  \|_{H^{\ell-1}}\big)\cr
&+CT^{\theta_{\ell}} \big( \|\psi\|_{X_S^{\ell+\frac 12,b}}+
 \|n\|_{X_{W_+}^{\ell,b}}+\|n_{t}\|_{X_{W_+}^{\ell-1,b}} \big)^2\label{apriori2},
\end{align}
with a power of $\theta_\ell$ as stated in Theorem \ref{thmrate}.

{\it Step2: Contradiction argument and lower bound.}

Let us  explain the contradiction argument 
 for a general evolution PDE with a quadratic nonlinearity, and an initial data $u_0$ belonging to some space $H$.  
Suppose  an a priori estimate of the form
\be\label{apriori}
\norm{u}_{X_{T}}\leq C\norm{u_{0}}_{H}+CT^{\theta}\norm{u}^{2}_{X_{T}},\ \ \theta >0
\ee
holds, 
where $\norm{.}_{X_T}$ is some appropriate  space-time norm. 
(This is   the a priori estimate \eqref{apriori2} with  $u=(\psi,n,n_t)$ and $H=\mathcal{H}_\ell$). 
Let $$\mathcal X(T,M)=\{ u: \ \ u(0)=u_{0},\ \ \norm{u}_{X_{T}}\leq M \}.$$  When performing an iteration argument in $\mathcal X(T,M)$, we would like to show
\begin{equation}\label{eM}
C\norm{u_{0}}_{H}+CT^{\theta}M^{2}\leq M
\end{equation}
to keep the iterates in $\mathcal X(T,M)$.   Local wellposedness  follows if \eqref{eM} holds (with, for example,
   $M=2C\norm{u_{0}}_{H}$, and  $T$ small enough so that $2CT^{\theta}M<1$).     
The relation between the spaces $X_T$ and $H$ is that $X_T$ must be imbedded in $C([-T,T],H)$ meaning that if $u$ belongs
to $X_T$, it must be a continuous function of $t \in [-T,T]$ with values in $H$.

Let $t^*$ be the maximal time of existence of solutions, that is 
 $$t^{*}=\sup \{T: \norm{u}_{X_{T}}<\infty\}.$$
 The blowup hypothesis implies that $t^*$ is finite.
Returning to  \eqref{eM}, let $0<t<t^*$ and consider  $u(t)$ as an initial condition. The following statement must hold:\\
If there exists some $M>0$ such that 
$ C\norm{u(t)}_{H}+C(T-t)^{\theta}M^{2}\leq M,$
then $T<t^{*}$.  Or equivalently: If  $T\geq t^{*}$, in particular $T=t^{*}$, then for all $M>0$
\[
C\norm{u(t)}_{H}+C(t^{*}-t)^{\theta}M^{2}> M.
\]
We now choose  $M=2C\norm{u(t)}_{H}$, then $\frac 12M+C(t^{*}-t)^{\theta}M^{2}> M $ or
\be
 C(t^{*}-t)^{\theta}M^{2}> \frac{M}{2} 
 \ee
 or equivalently 
 \be \label{blowuprate}
  \norm{u(t)}_{H}> c (t^{*}-t)^{-\theta}.
\ee
Hence, since we cannot continue the time of existence past time $t^\ast$, we have  a lower bound for the blowup rate of the norm $H$ as given by \eqref{blowuprate}.  

 Note that the conclusion is about the rate of blowup of the norm $H$ even though the iteration is performed in another norm.  One just needs the other norm to embed into $C([-T,T],H)$.

\section{The Vectorial Zakharov system}
There is no rigorous analysis of blowing up solutions of the full vectorial Zakharov system \eqref{Z1}-\eqref{Z3}
although wave collapse is expected when the initial conditions are large enough
 on the basis of numerical simulations and heuristic  arguments (see \cite{SS99} for review). Here we extend
 the results  of Section 2.3.3.

Recalling that for a vector valued function $\bf E$
\[
\Delta \bf E=\nabla (\nabla \cdot \bf E)-\nabla\times(\nabla \times \bf E),
\]
we can write \eqref{Z1}-\eqref{Z3} as
  \begin{align}
& i\partial_t {\bf E} + \alpha \Delta {\bf E}  +(1-\alpha)  
\nabla(\nabla \cdot {\bf E}) =n {\bf E}, \label{Z1a}  \\
&\partial_{tt} n - \Delta n = \Delta |{\bf E}|^2. \label{Z2a}  
\end{align}
The symbol of the Laplacian  is $\abs{\xi}^2$.  This together with the one time derivative determines the $\ip{ \tau +\abs{\xi}^2}$ weight in the $X^{s,b}$ norm for the NLS equation of the scalar Zakharov system.  To determine the weight that should appear in the $X^{s,b}$ norm for the NLS equation of the vectorial Zakharov system, one needs to determine the symbol of the spatial  linear operator appearing in the lhs of \eqref{Z1a}.  A simple calculation  leads to the matrices $M_2$ in two dimensions and $M_3$ in three dimensions given by 
\begin{align*}
M_2&=(1-\alpha)\left( \begin{array}{cc}
\xi_1^2 & \xi_1\xi_2\\
\xi_1\xi_2 &\xi_2^2\end{array} \right)+\alpha \abs{\xi}^2I_{2\times 2},\\
\intertext{and}
M_3&=(1-\alpha)\left( \begin{array}{ccc}
\xi_1^2 & \xi_1\xi_2  & \xi_1\xi_3\\
\xi_1\xi_2 & \xi_2^2   & \xi_2\xi_3\\
\xi_1\xi_3 & \xi_2\xi_3   & \xi_3^2 \end{array} \right)+\alpha \abs{\xi}^2I_{3\times 3}.
\end{align*}
where $I_{d\times d}$ is the $d\times d$ unit matrix.
It was observed by Tzvetkov \cite{T00} that the symbol of the operator given by matrix $M_d$ is actually equivalent to the symbol of the Laplacian.
\begin{lemma}\cite[Proposition 1]{T00}  Let $d=2,3$.  Then there exists a constant $C$ such that
\[
\abs{\xi}^2I_{d\times d}\leq M_d\leq C\abs{\xi}^2I_{d\times d}.
\]
 
\end{lemma}
Using this lemma, one can obtain a local well-posedness result for \eqref{Z1a}-\eqref{Z2a} \cite{T00} analogous to the scalar case, and a lower bound for the rate of blowup of Sobolev norms.

\begin{theorem}  Let $d=2,3$, and the  initial data $(E(0), n(0), n_t(0))$
 be  in   $ {\mathcal H_\ell} := (H^{\ell+1/2}(\R^d))^d  \times H^{\ell}(\R^d) \times H^{\ell-1}(\R^d)$,  
  $0\leq \ell \leq \frac d2 -  \frac 12 $.
  Assume that the solution $(E,n,n_t)$ blows up in a finite time $t^\ast <\infty$.
Then
\begin{equation}
\|{\mathbf E} (t)\|_{H^{\ell+1/2}} +\|n(t)\|_{H^{\ell}} + \|n_t(t)\|_{H^{\ell-1}}
> C(t^*-t)^{-\theta_\ell} \label{BUPrateVect}
\end{equation}
with $\theta_\ell =\frac14(4-d+2\ell)^-$, in dimension  $d=2 $ or  $d=3$.

\end{theorem}

Finally, as for the NLS equation,  it is of  interest to consider the influence of  additional dispersive terms
and their effect on blowing up solutions.  
In \cite{HS09}, Haas and Schukla consider the system
\begin{align}
& i\partial_t {\mathbf E} - \alpha\mbox{\boldmath $\nabla$} \times 
( \mbox{\boldmath $\nabla$}\times {\mathbf E})  +  
\mbox{\boldmath $\nabla$}(\mbox{\boldmath $\nabla$} \cdot {\mathbf E})   = n {\mathbf E} + 
\Gamma  \, \mbox{\boldmath $\nabla$}\Delta (\mbox{\boldmath $\nabla$} \cdot {\mathbf E})  \label{ZV1}\\
&\partial_{tt} n - \Delta  n = \Delta |{\mathbf E}|^2 -\Gamma \Delta^2 n. \label{ZV2}
\end{align}
which takes into account quantum corrections. 
The coefficient $\Gamma>0$ is  assumed  to be very small. 
In \cite{SSS09}, it is shown rigorously that quantum terms  arrest collapse    
in two and three dimensions,  for arbitrarily small values of the  parameter $\Gamma$. 

\begin{acknowledgment}
MC  is partially supported by grant  \#246255 from the Simons Foundation. CS  is partially supported by
NSERC through grant number 46179--13 and  Simons Foundation Fellowship
\#265059. 
\end{acknowledgment}

\end{document}